\documentclass[11pt,reqno,letterpaper]{amsart}
\usepackage{color}
\usepackage[colorlinks=true, allcolors=blue,backref=page]{hyperref}
\usepackage{amsmath, amssymb, amsthm}
\usepackage{mathrsfs}
\usepackage{mathtools}
\usepackage[noabbrev,capitalize,nameinlink]{cleveref}
\crefname{equation}{}{}
\usepackage{fullpage}
\usepackage[noadjust]{cite}
\usepackage{graphics}
\usepackage{pifont}
\usepackage{tikz}
\usepackage{bbm}
\usepackage[T1]{fontenc}

\usetikzlibrary{arrows.meta}

\usepackage{environ}
\usepackage{framed}
\usepackage{url}
\usepackage[linesnumbered,ruled,vlined]{algorithm2e}
\usepackage[noend]{algpseudocode}
\usepackage[labelfont=bf]{caption}
\usepackage{cite}
\usepackage{framed}
\usepackage[framemethod=tikz]{mdframed}
\usepackage{appendix}
\usepackage{graphicx}
\usepackage[textsize=tiny]{todonotes}
\usepackage{tcolorbox}
\usepackage{enumerate}
\allowdisplaybreaks[1]
\usepackage{enumerate}
\usepackage{stmaryrd}
\usepackage[margin=1in]{geometry}

\usepackage[shortlabels]{enumitem}
\crefformat{enumi}{#2#1#3}
\crefrangeformat{enumi}{#3#1#4 to~#5#2#6}
\crefmultiformat{enumi}{#2#1#3}
{ and~#2#1#3}{, #2#1#3}{ and~#2#1#3}

\DeclareSymbolFont{symbolsC}{U}{pxsyc}{m}{n}
\SetSymbolFont{symbolsC}{bold}{U}{pxsyc}{bx}{n}
\DeclareFontSubstitution{U}{pxsyc}{m}{n}
\DeclareMathSymbol{\medcircle}{\mathbin}{symbolsC}{7}

\crefname{algocf}{Algorithm}{Algorithms}

\crefname{equation}{}{} 
\AtBeginEnvironment{appendices}{\crefalias{section}{appendix}} 

\usepackage[color,final]{showkeys} 

\colorlet{refkey}{orange!20}
\colorlet{labelkey}{blue!30}

\crefname{algocf}{Algorithm}{Algorithms}

\numberwithin{equation}{section}
\newtheorem{theorem}{Theorem}[section]

\newtheorem{lemma}[theorem]{Lemma}

\crefname{claim}{Claim}{Claims}

\newtheorem*{question*}{Question}

\theoremstyle{definition}
\newtheorem{definition}[theorem]{Definition}

\newtheorem*{definition*}{Definition}

\theoremstyle{remark}
\newtheorem*{remark}{Remark}

\newtheorem{assumption}[theorem]{Assumption}

\newcommand{\mb}{\mathbb}
\newcommand{\mbf}{\mathbf}

\newcommand{\ol}{\overline}

\newcommand{\imod}[1]{~\mathrm{mod}~#1}

\newcommand{\eps}{\varepsilon}

\let\originalleft\left
\let\originalright\right
\renewcommand{\left}{\mathopen{}\mathclose\bgroup\originalleft}
\renewcommand{\right}{\aftergroup\egroup\originalright}

\allowdisplaybreaks

\newcommand{\ignore}[1]{}

\title{An explicit economical additive basis}

\author[A1]{Vishesh Jain}
\address{Department of Mathematics, Statistics, and Computer Science, University of Illinois Chicago, Chicago, IL, 60607 USA}
\email{visheshj@uic.edu}

\author[A2]{Huy Tuan Pham}
\address{Department of Mathematics, Stanford University,
Stanford, CA 94305, USA}
\email{huypham@stanford.edu}

\author[A3]{Mehtaab Sawhney}
\author[A4]{Dmitrii Zakharov}
\address{Department of Mathematics, Massachusetts Institute of Technology, Cambridge, MA 02139, USA}
\email{\{msawhney,zakhdm\}@mit.edu}

\begin{document}

\begin{abstract}
We present an explicit subset $A\subseteq \mb{N} = \{0,1,\ldots\}$ such that $A + A = \mb{N}$ and for all $\eps > 0$, 
\[\lim_{N\to \infty}\frac{\big|\big\{(n_1,n_2): n_1 + n_2 = N, (n_1,n_2)\in A^2\big\}\big|}{N^{\eps}} = 0.\]
This answers a question of Erd\H{o}s. 
\end{abstract}

\maketitle

\section{Introduction}\label{sec:intro}

Sidon asked (see \cite{Erd54, Sid32}) whether there exists a set $A\subseteq \mb{N}$ such that $A + A = \mb{N}$ (i.e.~$A$ is an additive basis of order $2$) and for all $\eps > 0$,
\[\lim_{N\to \infty}\frac{\big|\big\{(n_1,n_2): n_1 + n_2 = N, (n_1,n_2)\in A^2\big\}\big|}{N^{\eps}} = 0.\]
Erd\H{o}s \cite{Erd54} answered Sidon's question by showing that there exists an additive basis of order $2$ which, in fact, satisfies the stronger bound 
\[\limsup_{N\to \infty}\frac{\big|\big\{(n_1,n_2): n_1 + n_2 = N, (n_1,n_2)\in A^2\big\}\big|}{\log N} <\infty.\]
It is a major open problem whether there exists an additive basis of order $2$ for which the factor of $\log{N}$ in the denominator can be replaced by an absolute constant; Erd\H{o}s and Tur\'{a}n \cite{ET41} famously conjectured that this is impossible. 

Erd\H{o}s's proof of the existence of $A$ is randomized; in modern notation, one simply includes the number $n$ in the set $A$ with probability proportional to $C(\log n)^{1/2}n^{-1/2}$. Kolountzakis \cite{Kol95} derandomized (a variation of) Erd\H{o}s's proof in the sense that one can deterministically generate the elements of $A\cap \{0,\ldots, N\}$ in time $N^{O(1)}$. We remark that a number of variants of the original result of Erd\H{o}s have been developed including results of Erd\H{o}s and Tetali \cite{ET90} which prove the analogous result for higher order additive bases and results of Vu \cite{Vu00} regarding economical versions of Waring's theorem.

The focus of this work is on ``explicit'' constructions. Erd\H{o}s several times \cite{EG80,Erd95,Erd97} asked for an explicit set $A$ which affirmatively answers Sidon's question and, in fact, offered $\$100$ for a solution \cite{Erd97}. We note that if one takes $A$ to be the set of squares, then $A+A$ contains all primes which are $1\imod 4$ and by the divisor bound $A+A$ has multiplicities bounded by $N^{o(1)}$. Therefore, if one is willing to assume strong number theoretic conjectures, one can take $A$ to be the set of numbers $n$ which are within $O((\log n)^{O(1)})$ a square. The purpose of this note is to present an explicit construction 
unconditionally. 

Given a set $A$ which is either a subset of $\mb N$ or $\mb Z/q\mb Z$ for some $q$, we denote by $\sigma_A(n)$ the number of representations $n = a+a'$ or $n \equiv a+a' \imod q$, where $a, a' \in A$. 

\begin{theorem}\label{thm:main}
There is an explicit set $A \subset \mb N$ and absolute constants $C, c > 0$ such that for every $n \in \mb N$, we have $1 \le \sigma_A(n) \le C n^{c/\log \log n}$.
\end{theorem}

We briefly discuss the meaning of the word ``explicit''. By analogy with the long line of work on explicit Ramsey graphs, we adopt the convention that a construction is explicit if one may test membership $n\in A$ in $(\log n)^{O(1)}$-time i.e.~polynomial in the number of digits. In addition to satisfying this guarantee, our construction has the appealing feature that, given a suitable sequence of prime numbers, we can describe $A$ with a short explicit expression (see~\cref{eq:description-A} and the proof of \cref{lem:modular}).

We end this introduction with a brief overview of the (short) proof of \cref{thm:main}. A crucial ingredient in our work is a construction of Rusza \cite{Ruz90} which (for a prime $p \equiv 3,5 \imod 8$) produces a set $A_p\subseteq \mb{Z}/(p^2\mb{Z})$ such that $A_p + A_p = \mb{Z}/(p^2 \mb{Z})$ and $\sigma_{A_p}(r) = O(1)$ for all $r \in \mb{Z}/(p^2 \mb{Z})$. Given this, the key point is to consider a sequence of such primes $p_1<p_2<\ldots $ and define the set $A$ in terms of its expansion with respect to the generalized base $(p_1^2,p_2^2,\dots)$ essentially by forcing the $i^{th}$ digit in the expansion to belong to $A_{p_i}$. By working upward from the smallest digit, one can use the property that $A_{p} + A_p$ covers all residues in $\mb{Z}/(p^2 \mb{Z})$ to see that all natural numbers are represented. The fact that no number is represented too many times is similarly derived using that $A_p + A_p$ is ``flat'', in particular that the multiplicities are bounded by $p^{o(1)}$. We remark that generalized bases have been utilized in a variety of questions related to Sidon sets, including works of Ruzsa \cite{Ruz98}, Cilleruelo, Kiss, Ruzsa and Vinuesa \cite{CKRV10}, and Pillate \cite{Pil23}. 

We note that for the purpose of simply obtaining an upper bound of $N^{o(1)}$ in \cref{thm:main}, one can actually choose the primes $p$ sufficiently small (e.g. of size $\log\log\log N$ say) and find a suitable set $A_p$ by brute force enumeration. However Rusza's \cite{Ruz90} construction is ``strongly explicit'' (i.e.~membership can be tested in time $O((\log p)^{O(1)})$), which allows us to take larger primes and thereby obtain a better upper bound. 

The computational bottleneck (given Rusza's construction) is finding the smallest prime in an interval $[N,2N]$. Under strong number-theoretic conjectures (e.g. Cramer's conjecture) finding such a prime would take time $O((\log N)^{O(1)})$ due to the AKS primality testing algorithm \cite{AKS04}. Assuming this to be the case, we can choose the primes more carefully to obtain an improved upper bound of $\exp(O((\log N)^{1/2}))$ (see \cref{sec:cramer}). The limiting feature of our construction now is that in the top block, one is forced to allow ``all possibilities''.
We believe obtaining an explicit construction achieving $\sigma_A(N) \leq \exp((\log N)^{\eps})$ or better would be interesting. 

\subsection*{Notation}
Throughout this paper we let $[N] = \{0,\ldots, N-1\}$ and $\mb{N} = \{0,1,2,\ldots\}$. We let $\lfloor x\rfloor$ denote the largest integer less than or equal to $x$. We use standard asymptotic notation, e.g. $f \lesssim g$ if $|f(n)| \le C |g(n)|$ for some constant $C$ and all large enough $n$. We usually denote by $c, C$ absolute constants which may change from line to line.

\section{Proof of \cref{thm:main}}

We will require the notion of a generalized base. 
\begin{definition}\label{def:gen-base}
Let $\mbf{b} = (b_1,b_2,\ldots)$ be an infinite set of integers such that $b_i\ge 2$. Given any integer $x\in \mb{N}$ there exists a representation $x = \ol{a_n\ldots a_1}^{\mbf{b}}$ with $0\le a_i\le b_i - 1$ such that 
\[x = \sum_{i=1}^{n}a_i \prod_{j<i}b_j.\]
Here an empty product (when $i=1$) is treated as $1$.
\end{definition}
\begin{remark}
If one requires $a_n \neq 0$ (e.g. does not have leading zeros) the representation is unique. When $b_j = g$ for all $j$, we recover precisely the base-$g$ expansion.
\end{remark}

A crucial piece in our construction is an ``economical'' modular additive basis of order $2$ over $\mb{Z}/(p^2\mb{Z})$ due to Ruzsa \cite[Theorem~1]{Ruz90}; the precise constant $M$ in the result below has been studied in \cite{Che08}.  
\begin{lemma}\label{lem:modular}
There exists an absolute constant $M\ge 1$ such that the following holds. Consider a prime $p$ such that $p \equiv 3,5 \imod 8$. There exists a set $A_p\subseteq \mb{Z}/(p^2 \mb{Z})$ such that for all $r\in \mb{Z}/(p^2\mb{Z})$, we have $1\le \sigma_{A_p}(r)\le M$.
Furthermore, given $p$ and $x\in \mb{Z}/(p^2 \mb{Z})$, one can check whether $x\in A_p$ in time $O((\log p)^{O(1)})$.
\end{lemma}

For completeness, and especially in order to discuss the second part of the statement, we present the proof of \cref{lem:modular} in \cref{sec:mod-construc}.

We next need the following basic fact about deterministically finding primes, which is immediate via (say) the Sieve of Eratosthenes. Using a more sophisticated algorithm of Lagarias and Odlyzko \cite{LO87}, one may obtain a run time of $O(N^{1/2+o(1)})$ in the statement below; runtimes of the form $O(N^{o(1)})$ remain a major open problem. 
\begin{lemma}\label{lem:prime-prod}
Let $N \ge C_{\ref{lem:prime-prod}}$. Then one may produce the smallest prime $p\in [N,2N]$ such that $p \equiv 3\imod 8$ in time $O(N^{1+o(1)})$. 
\end{lemma}

We now are in position to give the proof of \cref{thm:main}.
\begin{proof}[{Proof of \cref{thm:main}}]
Let $f: \mb{N} \rightarrow \mb{N}$ be an arbitrary monotone increasing function such that $f(k) \ge C_0$ for some large constant $C_0$ and all $k$. Let $p_1 < p_2 < \ldots $ be a sequence of primes such that $p_k \equiv 3 \imod 8$ and $p_k \in [f(k), 2f(k))$ is the least such prime.\footnote{That such a prime exists for $f(k)$ larger than an absolute constant follows via the Siegel-Walfisz theorem (see e.g. \cite[Theorem~12.1]{Kuo19}).} Define ${\mbf b} = (b_1, b_2, \ldots)$ by setting $b_{k} = p_k^2$ for all $k \ge 1$. We are going to define our set $A$ in terms of its expansion in  the generalized base ${\mbf b}$. Namely, for each $k\ge 1$ let $A_k \subset \{0, 1, \ldots, p_k^2-1\}$ be the set from \cref{lem:modular} (where we lift elements $\imod p_k^2$ to their integer representatives) and consider the set 
\begin{equation}
\label{eq:description-A}
A = \bigcup_{k\ge 1}\{ \overline{a_k a_{k-1} \ldots a_1}^{\mbf{b}},  ~ a_j \in A_j\text{ for }j =1, \ldots, k-1, \text{ and }a_k \in \{0, \ldots, b_k-1\}\}.
\end{equation}

We begin by showing that $A + A = \mb{N}$. For any $n \in \mb{N}$, we (uniquely) write $n = \ol{n_k \ldots n_1}^{\mbf b}$ for some $k \geq 1$; we will construct the representation $n = a+a'$ for $a, a' \in A$ digit by digit. First, since $A_1$ is an order $2$ additive basis mod $b_1$, there exist $a_1, a_1'\in A_1$ such that $n_1 \equiv a_1 +a'_1 \imod b_1$. Let $c_1 = \big\lfloor\frac{a_1+a'_1}{b_1}\big\rfloor \in \{0,1\}$ be the carry bit. Next, there exist $a_2, a'_2\in A_2$ such that $n_2 - c_1 \equiv a_2 + a_2' \imod b_2$. As before, define the carry bit $c_2$ and continue in the same fashion to produce sequences of digits $a_1, \ldots, a_{k-1}$ and $a'_1, \ldots, a'_{k-1}$ and a carry bit $c_{k-1} \in \{0,1\}$. Finally, let $a_k = n_k - c_k \in \{0, \ldots, b_k-1\}$ and consider the elements
\[
a = \ol{a_k\ldots a_1}^{\mbf b}, ~~ a' = \ol{a'_{k-1} \ldots a'_1}^{\mbf b}.\]
By construction, we have $n = a+a'$ and $a, a' \in A$. 

Next, we bound the number of possible representations $n = a + a'$ with $a,a' \in A$. Write $n = \ol{n_k\dots n_1}^{\mbf b}$, $a = \ol{a_{\ell}\ldots a_1}^{\mbf b}$, and $a' = \ol{a'_{\ell'}\ldots a'_1}^{\mbf b}$, where $\ell, \ell' \le k$ are the digit lengths of $a$ and $a'$. We may assume that $\ell \le \ell'$ (this costs us a factor of 2 in the number of representations). Since $a, a' \in A$, we have $a_i \in A_i$ for $i \le \ell-1$ and $a'_i\in A_i$ for $i \le \ell'-1$ but the top digits $a_\ell$ (respectively, $a'_{\ell'}$) can be arbitrary elements of $\{0,\dots, b_{\ell}-1\}$ (respectively, $\{0,\dots, b_{\ell'}-1\}$). By \cref{lem:modular}, we can choose $a_1, a_1'$ such that $n_1 \equiv a_1+a_1' \imod b_1$ in at most $M$ ways. Given a choice of $a_1, a_1'$, there are at most $M$ pairs $a_2, a_2'$ with $n_2-c_1 \equiv a_2+a_2' \imod b_2$, where $c_1 = \big\lfloor\frac{a_1+a_1'}{b_1}\big\rfloor\in \{0,1\}$ is the carry. Continuing in this fashion for $j=1, \ldots, \ell-1$, we get that there are at most $M^{\ell-1}$ ways to fix the first $\ell-1$ digits $a_1,\dots, a_{\ell-1}$ and $a_1',\dots, a_{\ell-1}'$. We can fix $a_{\ell}$ and $a'_{\ell}$ in at most $b_\ell$ ways. Given this choice, the digits $a'_{\ell+1}, \ldots, a'_{\ell'}$ are uniquely determined. Putting this together, we obtain the following upper bound on the number of representations $n=a+a'$:
\begin{equation}\label{eq:1}
\sigma_A(n) \le 2 \sum_{\ell=1}^{k} b_\ell M^{\ell-1} \le 2 b_k M^k \le 8 f(k)^2 M^{k},
\end{equation}
where we used $b_k = p_k^2 \le (2f(k))^2$ and $M\ge 2$.
On the the other hand, $b_j = p_{j}^2\ge f(j)^2$ and so 
\begin{equation}\label{eq:2}
n \ge b_1 \ldots b_{k-1} \ge b_{\lfloor k/2\rfloor} \ldots b_{k-1} \ge f(\lfloor k/2\rfloor)^{k/2}.
\end{equation}
Hence, $k \le \frac{2\log n}{\log f(\lfloor k/2\rfloor)}$ and substituting this in \cref{eq:1}, we obtain the bound
\[\sigma_A(n) \le 4 f(k)^2 n^{c / \log f(\lfloor k/2\rfloor)}.\]
Note that the right hand side is $n^{o(1)}$ for any sufficiently slowly growing function $f$.
Owing to the computational considerations in the next paragraph, we take $f(k) = k$ which leads to $k \lesssim \frac{\log n}{\log \log n}$ and $\sigma_A(n) \lesssim n^{c/\log \log n}$.

Finally, we quickly verify that testing membership  $a \in A$ can be done in time $O((\log a)^{O(1)})$. Indeed, given $a \in \mb N$, we can compute all primes $p_k$ for $k \le c\log a$ in time $(\log a)^{O(1)}$ (\cref{lem:prime-prod}), compute the base $\mbf b$ expansion $a = \ol{a_k\ldots a_1}$ in time $(\log a)^{O(1)}$, and check that $a_j \in A_j$ for $j=1, \ldots, k-1$ in time $O(k (\log f(k))^{O(1)})$ using \cref{lem:modular}. 
\end{proof}

\subsection{Modular construction and computational details}\label{sec:mod-construc}
We record the proof of \cref{lem:modular}, following Ruzsa \cite{Ruz90}. 
For $n \in \mb Z, p \in \mb N$ we write $(n \imod p)$ for the unique $n' \in \{0, 1, \ldots, p-1\}$ congruent to $n$ modulo $p$.
The following is exactly \cite[Lemma~3.1]{Ruz90}. 
\begin{lemma}\label{lem:basic-cont}
Let $p\equiv 3,5\imod 8$. Define $B_p \subseteq \{0,\dots, 2p^2\}$ by
\begin{align*}
B_p &= \{x + 2p (3x^2 \imod p): x\in \{0,\ldots,(p-1)\}\}\\
&\cup \{x + 2p (4x^2 \imod p): x\in \{0,\ldots,(p-1)\}\} \\
&\cup \{x + 2p (6x^2 \imod p): x\in \{0,\ldots,(p-1)\}\}. 
\end{align*}
We have $\sup_{n\in \mb{Z}}\sigma_{B_p}(n) \le 18$ and furthermore, for all $0\le n<p^2$, at least one of the six numbers
\[n-p, n, n+p, n+p^2-p, n+p^2,n+p^2+p\]
appears in the set $B_p + B_p$.
\end{lemma}

Given \cref{lem:basic-cont}, we prove \cref{lem:modular}.
\begin{proof}[{Proof of \cref{lem:modular}}]
Let $B_p' = B_p + \{-p, 0, p\}$ (viewed as a subset of $\mb{Z}$) and set 
\[A_p = \{x \imod p^2: x\in B_p'\} \subseteq \mb{Z}/(p^2 \mb{Z}).\]
Applying \cref{lem:basic-cont}, we immediately have:
\begin{itemize}
    \item $B_p' + B_p' \subseteq [-2p,5p^2]$
    \item For all $0\le n<p^2$, one of $n$ or $n+p^2$ appears in $B_p' + B_p'$
    \item We have that $\sup_{n\in \mb{Z}}\sigma_{B_p'}(n) \le 9 \cdot 18 = 162$. 
\end{itemize}
Noting that $n\equiv n + p^2 \imod p^2$, it follows that $A_p + A_p = \mb{Z}/(p^2\mb{Z})$. Furthermore we have that $\sup_{n\in \mb{Z}}\sigma_{A_p}(n) \le 6 \cdot 9 \cdot 18 = 594$; this is immediate as $B_p' + B_p' \subseteq [-2p,5p^2]$ and there are at most $6$ representatives in this interval for a given residue modulo $p^2$. 

We now discuss the time complexity of testing membership in $A_p$. Given $p$, and $x\in \mb{Z}/(p^2 \mb{Z})$, we consider the unique representative $x' \in \{0,\ldots, p^2-1\}$. Noting that $B_p'\subseteq [-p, 2p^2 + p]$, it suffices by construction to test whether at least one of $x' - p^2, x', x'+p^2, x'+2p^2$ is in $B_p'$. This is equivalent to checking whether at least one of $x' + \{-p^2,0,p^2,2p^2\} + \{-p,0,p\}$ is in $B_p$; in particular, one of at most $12$ distinct given elements is in $B_p$. 

To test whether $y\in \{0,\ldots,2p^2\}$ is in $B_p$ amounts to testing whether $y = z + 2p (3z^2 \imod p)$, $y = z + 2p (4z^2 \imod p)$, or $y = z + 2p (6z^2 \imod p)$ for an integer $z\in \{0,\ldots, p-1\}$. Given $y$, the ``candidate'' $z$ is forced to be the unique number in $\{0,\dots,p-1\}$ equivalent to $y \imod p$ and we can then simply compute $(3z^2 \imod p)$, $(4z^2 \imod p)$, and $(6z^2 \imod p)$. This procedure clearly takes time $O((\log p)^{O(1)})$. 
\end{proof}

\subsection{Assuming deterministic polynomial time algorithms for locating primes}\label{sec:cramer}
For the remainder of this section we will operate under the following assumption.
\begin{assumption}\label{asp:prime}
There exists a deterministic algorithm which outputs the least prime which is $3 \imod 8$ in the interval $[N,2N]$ in time $O((\log N)^{O(1)})$.
\end{assumption}

To obtain a better upper bound on $\sigma_A(n)$, we take the function $f$ in the proof of \cref{thm:main} to be $f(k) = \exp(c k)$. It follows from (\ref{eq:1}) and (\ref{eq:2}) that $\sigma_A(n) \lesssim \exp(Ck)$ and $n \gtrsim \exp(c k^2)$, thus giving the bound $\sigma_A(n) \lesssim \exp(C \sqrt{n})$. To test membership, we need to construct primes $p$ of order at most $\exp(ck) \approx \exp(c\sqrt{\log n})$ which can be done in $(\log n)^{O(1)}$-time under \cref{asp:prime}.

\section*{Acknowledgements}
V.J.~is supported by NSF CAREER award
DMS-2237646. H.P.~is supported by a Clay Research Fellowship and a Stanford Science Fellowship. M.S.~is supported by NSF Graduate Research Fellowship Program DGE-2141064. D.Z.~is supported by the Jane Street Graduate Fellowship. We thank Zach Hunter and S\'andor Kiss for carefully reading the manuscript and suggesting improvements and references.

\bibliographystyle{amsplain0.bst}
\bibliography{main.bib}

\end{document}